\documentclass[12pt]{article}
\usepackage{amsmath}

\usepackage[psamsfonts]{amsfonts}

\newcommand{\be}{\begin{equation}} 
\newcommand{\bea}{\begin{eqnarray}}
\newcommand{\ee}{\end{equation}}
\newcommand{\beas}{\begin{eqnarray*}}
\newcommand{\eea}{\end{eqnarray}}
\newcommand{\eeas}{\end{eqnarray*}}

\def\Z{{\bf Z}}
\def\C{{\bf C}}

\newtheorem{formula}{Theorem}[section]

\newtheorem{corollary}[formula]{Corollary}

\begin{document}

\newpage
\renewcommand{\baselinestretch}{1.3}
\normalsize

\Large
\begin{center}

{\bf{Counting Conjugacy Classes of Elements of Finite Order in Lie Groups}}

\large

\vskip 1cm 
Tamar Friedmann$^*$ and Richard P. Stanley$^\dagger$
\vskip .6cm 
\normalsize
{\it $^*$University of Rochester, Rochester, NY}

{\it $^\dagger$Massachusetts Institute of Technology, Cambridge, MA}

\vskip .5cm 
\normalsize
\abstract{
Using combinatorial techniques, we answer two questions about simple classical Lie groups. Define $N(G,m)$ to be the number of conjugacy classes of elements of finite order 
$m$ in a Lie group $G$, and $N(G,m,s)$ to be the number of such classes whose elements have $s$ distinct eigenvalues or 
conjugate pairs of eigenvalues. What is $N(G,m)$ for $G$ a unitary, orthogonal, or symplectic group? What is $N(G,m,s)$ for these groups? For some cases, the first question was answered a few decades ago via group-theoretic techniques. It appears that the second question has not been asked before; here it is inspired by questions related to enumeration of vacua in string theory. Our combinatorial methods allow us to answer both questions.
}
\end{center}
\normalsize
\vskip .3cm



\vskip .3cm

\section{Introduction}

Given a group $G$ of linear transformations and integers $m$ and $s$,
let 
\be E(G,m)=\{ x\in G \; |\;  x^m=1\}. \ee
Also let
\be \label{egms1} E(G,m,s)=  \{ x\in E(G,m)\;  | \; x \mbox{ has } s \mbox{ distinct eigenvalues} \} \ee
for $G$ a unitary group, and
\be \label{egms2} E(G,m,s)= \{ x\in E(G,m) \; |\;  x \mbox{ has } s \mbox{ distinct conjugate pairs of eigenvalues}\} \ee
for $G$ a symplectic or orthogonal group,
and let 
\beas N(G,m)&=& \mbox{number of conjugacy classes of } G \mbox{ in } E(G,m), \\
N(G,m,s)&=& \mbox{number of conjugacy classes of } G \mbox{ in } E(G,m,s).
\eeas
For  $\Gamma$ any finitely generated abelian group and $G$  a Lie
group, one can consider the space of homomorphisms Hom$(\Gamma, G)$
and the space of representations of $\Gamma$ in $G$, 
that is, consider
\[ \mathrm{Rep}(\Gamma, G)\equiv \mathrm{Hom}(\Gamma, G)/G \]
(where $G$ acts by conjugation); using this
notation, 
\[ E(G,m)=\mathrm{Hom}(\Z/m\Z, G)\] 
and
\[N(G,m)=|\mathrm{Rep}(\Z/m\Z , G)|.\]  
For 
the case $\Gamma = \Z ^n$,  the spaces $\mathrm{Hom} (\Z^n, G)$ and
$\mathrm{Rep}(\Z^n, 
G)$ have been studied for various Lie groups $G$ in \cite{BFM, AC,
  ACG, BJS} (and references therein), where there has been interest in
their number of path-connected components and their cohomology groups.  

It is the purpose  of this paper to compute $N(G,m)=|\mathrm{Rep}(\Z/m\Z , G)|$
and $N(G,m,s)$ for $G$ a unitary, orthogonal, or symplectic group.
Unlike $\mathrm{Rep}(\Gamma, G)$ for $\Gamma = \Z^n$, the representation space
$\mathrm{Rep}(\Z/m\Z , G)$ is a finite set,  so we can  count its number of
elements. The results are summarized in Table 1.  

The numbers $N(G,m,s)$ have never been studied before in the
mathematical literature.  
What motivated their definition, as well as the definition of
$N(G,m)$, was the need to find a formula for the number of certain
vacua in the quantum moduli space of M-theory compactifications on
manifolds of $G_2$ holonomy. In that context, the numbers $N(SU(p),q)$
and $N(SU(p), q, s)$, where $q$ and $p$ are relatively prime, were
computed in \cite{Friedmann:2002ct}. These numbers are related to
symmetry breaking patterns in grand unified theories, with the number
$N(SU(p),q,s)$ being particularly significant as $s$ is related to the
number of massless fields 
in the gauge theory that remains after the symmetry breaking. The
connections with symmetry breaking patterns arise from the fact that
if $M$ is a manifold and $\pi_1 (M)$ is its fundamental group, then
$\mathrm{Rep}(\pi _1(M), G)$ is the moduli space of isomorphism classes of flat
connections on principal $G$-bundles over $M$; in grand unified
theories arising from string or M-theory, these flat connections
(called Wilson lines) serve as a symmetry breaking mechanism. For more
on the physical applications and implications of these numbers, see
\cite{Friedmann1}. 

As for $N(G,m)$, certain cases have been studied previously in the
mathematical literature, using different techniques than ours. Two of
the quantities we derive, Theorems~\ref{qinsup} and \ref{qinspn}, were
obtained in \cite{81h:20052,86h:22010} using the full machinery of Lie
structure theory with a generating function approach; in
\cite{pianzola1, pianzola2}, the case of certain prime power orders is
computed; and in \cite{lossers}, Theorem~\ref{minsun} is
obtained. Our methods are different; they are purely combinatorial and
direct, and apply not only to simply connected or adjoint groups as in
\cite{81h:20052, 86h:22010}, so we are able to derive formulas for
$O(n)$, $SO(n)$, and $U(n)$ alongside those for $SU(n)$ and $Sp(n)$.

Other aspects of elements of finite order in Lie groups have been
studied. See for example \cite{other3, other1, other2, other4,
  other5}.

In addition to the quantities $N(G,m)$ and $N(G,m,s)$, which count
conjugacy classes of elements of any order dividing $m$, we  consider
also conjugacy classes of elements of exact order $m$ in $G$: let 

$$ F(G,m)=\{ x\in G \; |\;  x^m=1, x^n\neq 1 \mbox{ for all } n<m \}. $$
Also let
$$ F(G,m,s)= \{ x\in F(G,m) \; | \; x \mbox{ has } s \mbox{ distinct eigenvalues}\} $$
for $G$ a unitary group, and
$$ F(G,m,s)= \{ x\in F(G,m) \; | \; x \mbox{ has } s \mbox{ distinct conjugate pairs of eigenvalues} \}$$
for $G$ a symplectic or orthogonal group,
and let
\beas K(G,m)&=& \mbox{number of conjugacy classes of } G \mbox{ in } F(G,m), \\
K(G,m,s)&=& \mbox{number of conjugacy classes of } G \mbox{ in } F(G,m,s).
\eeas
Since
\beas N(G,m)&=&\sum _{d|m} K(G,d) ,\\
N(G,m,s)&=& \sum _{d|m}K(G,d,s) ,
\eeas
we have, by the M\"obius inversion formula, 
\bea \label{kgm} K(G,m)&=&\sum _{d|m} \mu (d) N(G,\frac{m}{d}) , \\ \label{kgms}
K(G,m,s)&=& \sum _{d|m} \mu (d) N(G,{m\over d},s) ,
\eea
where $\mu(d)$ is the M\"{o}bius function. 

The reader is invited to obtain $K(G,m)$ and $K(G,m,s)$ from Table 1
and equations (\ref{kgm}) and (\ref{kgms}) above.

\[ \renewcommand{\baselinestretch}{2.5}
\normalsize
\begin{array}{|cc| c| c|}
\hline
\multicolumn{4}{|c|}{\mbox{{\bf Table 1: Number of conjugacy classes
      of elements of finite order in Lie groups}}}\\ 
\hline
G& m & \hskip .5cm N(G,m)&\hskip .5cm N(G,m,s)\\
\hline
U(n)&\mbox{any}& {n+m-1\choose m-1}  &{s\over n} {n\choose s}{m\choose s} \\  
SU(n)& (n,m)=1& \frac{1}{m}{n+m-1\choose n}& {s \over nm} {n\choose
  s}{m\choose s} \\ 
& \mbox{any}&{1\over m}\sum\limits _{d|(n,m)}\phi(d) {(n+m-d)/ d
  \choose n/ d}&  {1\over m}\sum\limits _{d|(n,m)} \sum \limits_{
  j\geq 0}\phi(d) {(n+m-jd-d)/ d  \choose (n-jd)/ d}{m/d \choose
  j}{jd\choose s} (-1)^{j+s} \\ 
Sp(n)&\mbox{any}&   { n+[\frac{m}{2}]\choose n}  & {s\over n} { n
  \choose s} { [\frac{m}{2}]+1\choose s}\\  
SO(2n+1)&\mbox{any}& { n+[\frac{m}{2}] \choose n }& {s\over n} { n
  \choose s }  { \left [ {m\over 2}\right ] +1\choose s }\\ 
O(2n+1)&2k+1& { n+[\frac{m}{2}]\choose n } & {s\over n} { n \choose s
} { \left [ {m\over 2}\right ] +1\choose s }\\ 
O(2n)&2k+1& { n+[\frac{m}{2}] \choose n }& {s\over n} { n \choose s }{
  \left [ {m\over 2}\right ] +1\choose s} \\ 
SO(2n)&2k+1&  
{ n+[\frac{m}{2}]-1\choose n-1} \frac{n+m-1}{n} 
& 
{s\over n} { n \choose s } {\left [ {m\over 2}\right ] \choose s }
{m+1-s\over \left [ {m\over 2}\right ]+1-s } \\ 
O(2n+1)&2k&2{ n+\frac{m}{2}\choose n} & {2s\over n} { n \choose s }{
  {m\over 2}+1 \choose s}\\ 
O(2n)&2k&  
{ n+\frac{m}{2} -1\choose n-1 }{4n+m\over 2n}
& 
{2n-s-1\over n-s}{ n-2 \choose s -1 }{ {m\over 2}+1 \choose s}\\
SO(2n)&2k&  
{ n+\frac{m}{2}\choose n}+{ n+\frac{m}{2}-2\choose n }& {s\over n} { n
  \choose s }  \left [ 2{  {m\over 2}\choose s } + { {m\over
      2}-1\choose s-2 } \right ] 
\\
\hline
\end{array}\]

\renewcommand{\baselinestretch}{1.3}
\normalsize

\newpage
\section{Counting conjugacy classes in unitary groups}
We begin with $N(U(n), m)$, with no conditions on the integers $m$ and $n$. 
Since every element of $U(n)$ is  diagonalizable, every conjugacy
class has  diagonal elements. The diagonal entries are $m^{th}$ roots
of unity, $e^{2\pi i k_j/m}$, $k_j=0, \ldots, m-1$, and $j=1, \ldots ,
n$. In each conjugacy class there is a unique diagonal element for
which  the  diagonal entries are ordered so that the $k_j$  are
nondecreasing with $j$. Therefore, $N(U(n), m)$ is the number of such
diagonal matrices with nondecreasing $k_j$. 

Let 
$\{n_k\}=( n_0, \ldots , n_{m-1})$, $\sum_{k=0}^{m-1} n_k = n$ with
$n_k\geq 0$. Such a sequence is a weak $m$-composition of $n$, and it
is well-known that there are ${n+m-1 \choose m-1}$ such sequences
\cite{stanley}.  
There is a bijective map between such sequences and diagonal matrices
in $U(n)$ with ordered entries:  $\{n_k\}$ corresponds to the diagonal
$U(n)$ matrix with $n_k$ repetitions of the eigenvalue $e^{2\pi i
  k/m}$: 
\be \mbox{diag} ( \underbrace{1, 1, \cdots , 1}_{n_0},
\underbrace{e^{2\pi i/m}, \cdots , e^{2\pi i/m}}_{n_1}, \cdots ,
\underbrace{e^{2(m-1)\pi i/m}, \cdots , e^{2(m-1)\pi
    i/m}}_{n_{m-1}})~. 
\ee

 Thus $N(U(n), m)$ is the number of weak $m$-compositions of $n$, so
 we obtain the following formula.
\begin{formula}\label{qinup} For any positive integers $n$ and $m$, 
\be N(U(n),m)={n+m-1 \choose m-1} \ee
\end{formula}
Note that $N(U(n), m)$ is also the number of inequivalent unitary
representations of $\Z/m\Z$ of dimension $n$.  

Now we turn to the special unitary group $SU(p)$, and calculate
$N(SU(p),q)$ where $(p,q)=1$. Given a sequence $\{n_k\}$, $k=0, \ldots
q-1$ with $\sum_{k=0}^{q-1} n_k = p$, $n_k\geq 0$ (i.e. a weak
$q$-composition of $p$),  the determinant of the corresponding matrix
$x$  is $\exp {2\pi i \over q}\left (\sum_{k=0}^{q-1} kn_k \right )$,
so the condition $\det x=1$ requires  $\sum_k kn_k \equiv 0 \mbox{
  mod } q$. Thus for a weak $q$-composition of $p$ to determine a
matrix in $SU(p)$, we need  $\sum_k kn_k \equiv 0 \mbox{ mod } q$.  

We now show the family of weak $q$-compositions of $p$ are partitioned
into sets of size $q$ where in each such set there is exactly one such
composition with $\sum kn_k \equiv 0$. Consider the $q$  distinct
sequences 
\be \{n_k^{(j)}\}=\{n_{k+j}\} \hskip 1cm j=0, 1, \ldots , q-1~, \mbox{
  indices are understood mod q.} 
\ee
(The only way for the sequences not to be distinct is if all $n_k$
were equal, which would imply $qn_k=p$, impossible when $(p,q)=1$).  
The determinant of the matrix $x_j$  corresponding to the $j^{th}$
sequence is $\exp {2\pi i \over q}\left (\sum_{k=0}^{q-1}  kn_{k+j}
\right )$. Since $(p,q)=1$ and  
\be \sum_{k=0}^{q-1} kn_{k+j}-\sum_{k=0}^{q-1} kn_{k+j+1} \equiv p
\mbox{ mod } q~,\ee 
exactly one of the $q$ values of $j$ gives the sum $\sum_k
kn_{k+j}\equiv 0 \mbox{ mod } q$, so  $\det x_j =1$ for that value
of $j$. We therefore get the next result.
\begin{formula}\label{qinsup} For $(p,q)=1$, 
\be \label{qinsupder}
N(SU(p),q)=\frac{1}{q} {p+q-1 \choose q-1}= \frac{(p+q-1)!}{p!\, q!}~. 
\ee 
\end{formula}

Now we turn to counting conjugacy classes whose elements have a given
number $s$ of distinct eigenvalues. We begin with $N(U(n),m,s)$. A
$U(n)$ matrix with $s$ distinct eigenvalues (which has centralizer of
the form $\Pi _{i=1}^s U(n_i)$) corresponds to a sequence $\{n_a\} =
(n_1, \ldots n_s)$, $\sum _{a=1}^s n_a = n$, $n_a\geq 1$. Such a
sequence is an $s$-composition of $n$ and there are  ${n-1\choose
  s-1}$ such sequences \cite{stanley}). There are also 
${m\choose s}$
ways to choose the $s$ eigenvalues themselves. We therefore obtain the
following formula.

\begin{formula} \label{qsinup}For any positive integers $n$ and $m$,
$$ N(U(n),m,s)=\binom{n-1}{s-1}\binom ms = \frac{s}{n}
\binom ns\binom ms~.$$
\end{formula}
For the special unitary group, again we impose $(p,q)=1$. Given an
$s$-composition of $p$, $\{n_a\} = (n_1, \ldots n_s)$, $\sum _{a=1}^s
n_a = p$, $n_a>0$, consider  $\{ \lambda _a \} = (\lambda _1 , \ldots
, \lambda _s)$ where $\lambda _a \in \{0, \ldots , q-1 \}$ determine
the eigenvalues $e^{2\pi i\lambda _a\over q}$ with multiplicity $n_a$
of the corresponding matrix. Arrange the $\binom qs s! \hskip .1cm$
possibilities for  $\{ 
\lambda _a \} $ in sets of size $q$ given by 
\be \{ \lambda _a^{(j)}\}=(\lambda _1+j, \ldots , \lambda _s+j),\; \;
j=0,\ldots ,q-1 \hskip .5cm (\mbox{all numbers are understood mod
}q). \ee 
The determinant of the matrix $x_j$ corresponding to the $j^{th}$
choice is  
\[ \exp {2\pi i \over q}\left (\sum _{a=1}^s n_a (\lambda _a +j)\right
) .\] 
Since $(p,q)=1$ and 
\[ \sum _a n_a (\lambda _a+j)-\sum _a n_a(\lambda +j+1) = p,\]
exactly one of the $q$ matrices has determinant $1$. Since so far
neither the $\lambda _a$'s nor the $n_a$'s have been ordered, once we
arrange the eigenvalues to have  increasing $\lambda _a$'s, each
matrix would appear $s!$ times. Dividing by $s!q$, we obtain  the following formula.

\begin{formula} \label{qsinsup} For $(p,q)=1$,
\be N(SU(p),q,s)=\frac{1}{q} \binom{p-1}{s-1}\binom qs =
    \frac{s}{pq}\binom ps\binom qs ~. \ee 
\end{formula}

From Theorems \ref{qinsup} and \ref{qsinsup}, we deduce an intriguing
symmetry between $p$ and $q$. 
\begin{corollary} For $(p,q)=1$, 
\beas N(SU(p),q)&=&N(SU(q),p);\\ N(SU(p),q,s)&=&N(SU(q),p,s).
\eeas
\end{corollary}
This symmetry has implications involving dualities of gauge theories;
see \cite{Friedmann1}.  

It is clear that for any $G$ and $m$, we must have
\be \label{sum}\sum _s N(G,m,s)=N(G,m).\ee 
Since $N(G,m,s)=0$ when $s>m$, the sum is finite. 
Applying equation~(\ref{sum})  to $G=U(n)$ gives 
\be \label{id} \sum _s {n-1 \choose s-1}{m\choose s} = {n+m-1\choose m-1}, \ee
which  is a special case of the Chu-Vandermonde identity \cite{stanley}. 

We may also obtain both $N(SU(n),m)$ and $N(SU(n),m,s)$ without
requiring $(n,m)= 1$ via a generating function approach. Let 
\[ F(x,t,u)= \prod _{k=0}^{m-1} \left (1+u\sum_{a=1}^\infty (t^k x)^a \right ).
\]
A typical term in $F(x,t,u)$ is 
\[x^{\sum n_k}\, t^{\sum kn_k}\, u^s ,\] 
where $n_k$, $k=0, \ldots, m-1$ are nonnegative integers and $s$ is
the number of $k$'s for which $n_k\neq 0$. If $\sum n_k=n$ and $\sum
kn_k\equiv 0$ mod $m$ then the sequence $\{ n_k\}$ corresponds to a
diagonal $SU(n)$ matrix of order $m$ with $s$ distinct eigenvalues. To
pick out the terms in $F(x,t,u)$ for which $\sum kn_k\equiv 0$ mod
$m$, let $\zeta = \exp {2\pi i /m}$ and recall  
\[{1\over m}\sum _{j=0}^{m-1}\zeta ^{jb}=\left \{ \begin{array}{l} 1,\
  \mbox{ if } m|b  \\ 0,\ \mbox{ else} \end{array} \right . , 
\]
so 
$$ G(x,u)= {1\over m}\sum_{j=0}^{m-1}F(x,\zeta^j,u)=\sum_{n,s}N(SU(n),m,s)x^nu^s . $$
Rewriting  
\[ 1+u\sum_{a=1}^\infty (t^k x)^a = (1-u)+{u\over
  1-t^kx} = {1-t^k(1-u)x\over 1-t^kx},\] 
we have
$$ G(x,u)={1\over m}\sum_{j=0}^{m-1}\prod _{k=0}^{m-1}
    {1-\zeta^{kj}(1-u)x \over 1-\zeta^{kj}x} .$$
For $\zeta ^j$ a primitive $d^{th}$ root of unity, we have the
factorization $1-x^d=\prod_{l=0}^{d-1} (1-\zeta^{jl}x)$. Since
$\zeta ^j$, $j=0, \ldots, m-1$ is a primitive $d^{th}$ root of unity
$\phi (d)$ times, where $\phi(d)$ is Euler's  function, we have 
$$ G(x,u)={1\over m}\sum_{d|m}\phi(d) {\left [ 1-(1-u)^dx^d\right ]
  ^{m/d} \over (1-x^d)^{m/d} }. $$
Expanding in binomial series gives

\[ G(x,u)= {1\over m}\sum_{d|m}\phi(d) \sum _{k,j,l\geq 0} {k+m/d-1
  \choose k}{m/d \choose j}{jd\choose l}(-1)^{j+l}\, x^{d(k+j)}u^l .\] 

Setting  $d(k+j)=n$ and $l=s$ yields the next theorem.
\begin{formula}\label{msinsun} For any positive integers $n,m$, and $s$, 
$$ N(SU(n),m,s)= {1\over m}\sum_{d|(n,m)}\sum _{j\geq 0}\phi(d)
    {n/d+m/d - j -1 \choose n/d-j}{m/d\choose j}{jd\choose s}
    (-1)^{j+s}. 
$$
\end{formula}
We may deduce from Theorems~\ref{msinsun} and \ref{qsinsup} that for
$(p,q)=1$,  
$$ {1\over q}\sum _{j\geq 0} {p+q-j-1 \choose p-j}{q\choose
  j}{j\choose s}(-1)^{j+s}={s\over pq}{p\choose s}{q\choose s} .$$
For $N(SU(n),m)$ we apply equation~(\ref{sum}), or equivalently set
$u=1$ in $G(x,u)$, and obtain (see also \cite{lossers}) the next result.
\begin{formula}\label{minsun}  For any positive integers $n$ and $m$, 
$$ N(SU(n),m)={1\over m}\sum_{d|(n,m)}\phi(d) {n/d+m/d  -1 \choose
  n/d}. $$
\end{formula}

\section{Counting conjugacy classes in symplectic groups}
The diagonal elements of $U(n)$ and $SU(p)$ that we counted in the
previous section belong to the maximal tori of those groups. For
$\mathrm{Sp}(n)\equiv \mathrm{Sp}(n,\C)\cap U(2n)$, the maximal torus is  
\be \label{spntorus} T_{\mathrm{Sp}(n)}=\left \{ (e^{2\pi i\theta _1}, \ldots ,
e^{2\pi i\theta _n},  e^{-2\pi i\theta _1}, \ldots ,  e^{-2\pi i\theta
  _n}) \right \}. \ee 
Since $\mathrm{Sp}(n)$ is compact and connected, we have
$\mathrm{Sp}(n)=\bigcup_{x\in G} 
\;  xT_{\mathrm{Sp}(n)}x^{-1}$. Hence, every element $x\in G$ can be conjugated
into the torus, so 
every conjugacy class has elements in $T_{\mathrm{Sp}(n)}$. 
Any two elements $x$ and $x'$ of $T_{\mathrm{Sp}(n)}$ that differ only by
$\theta _l'=-\theta_l$ for some $l$'s are in the same conjugacy class;
the symplectic matrix $E_{l, n+l}-E_{n+l, l}$, where
$(E_{ab})_{cd}=\delta_{ac}\delta_{bd}$, conjugates them. So a
conjugacy class is fully determined by $n$ values of $\theta _l$
restricted to $[0,1/2]$.  

Conjugacy classes of elements of order $m$ have a unique element in
$T_{\mathrm{Sp}(n)}$ such that $\theta _l\in {1\over m}(0,1,\ldots , [{m\over
    2}])$ and the $\theta _l$ are nondecreasing as $i$ runs from 1 to
$n$. Following the arguments leading to Theorem \ref{qinup}, and
noting that here we have weak $([{m\over 2}]+1)$-compositions of $n$,
rather than weak $m$-compositions of $n$, we obtain our next theorem.

\begin{formula} \label{qinspn} For any positive integers $n$ and $m$,
$$ N(\mathrm{Sp}(n),m)=\binom{n+[\frac{m}{2}]}
{[\frac m2]}.$$
\end{formula}
We now consider $N(\mathrm{Sp}(n), m, s)$ where $s$ denotes the number of
complex conjugate pairs of eigenvalues.  Following the arguments
leading to Theorem  \ref{qsinup}, but replacing $m$ by $([{m\over
    2}]+1)$, we obtain the next result.

\begin{formula}\label{qsinspn} For any positive integers $n$, $m$, and $s$, 
$$ N(\mathrm{Sp}(n), m, s)=\binom{n-1}{s-1} \binom{[\frac{m}{2}]+1}{s}.$$
\end{formula}

\section{Counting conjugacy classes in orthogonal groups}
 
The maximal tori of the different orthogonal groups depend  on the
parity of $l$ in $SO(l)$ or $O(l)$ and also on whether the orthogonal
group is special or not: 
\bea
 \label{so2ntor}
T_{SO(2n)}&=&\left \{ \mbox{diag} (A(\theta _1), A(\theta _2), \ldots , A(\theta _n) )  \right \}, \\
\label{so2n+1tor}
T_{SO(2n+1)}&=&\left \{ \mbox{diag} (A(\theta _1), A(\theta _2), \ldots , A(\theta _n), 1 )  \right \}, \\
\label{o2ntor}
T_{O(2n)}&=&\left \{ \begin{array}{l}T_{1,even}=\mbox{diag} (A(\theta _1), A(\theta _2), \ldots , A(\theta _n) ) \\  T_{2,even}=\mbox{diag}(A(\theta _1), A(\theta _2), \ldots , A(\theta _{n-1}), B ) \end{array} \right \},\\
\label{o2n+1tor}
T_{O(2n+1)}&=&\left \{ \begin{array}{l} T_{1,\mathrm{odd}}= \mbox{diag}
  (A(\theta _1), A(\theta _2), \ldots , A(\theta _n),  1 )
  \\ T_{2,\mathrm{odd}}= \mbox{diag} (A(\theta _1), A(\theta _2), \ldots ,
  A(\theta _n),  -1 )\end{array}\right \}, 
\eea
where 
\be A(\theta)=\left ( \begin{array}{cc} \cos 2\pi \theta & \sin 2\pi
  \theta \\ -\sin 2\pi \theta & \cos 2\pi \theta  \end{array} \right
)\; ; \; B=\left ( \begin{array}{cc} 1 & 0 \\  0& -1  
\end{array} \right ).
\ee
In equations (\ref{o2ntor}) and (\ref{o2n+1tor}), the maximal torus is made of two parts. The first has elements of determinant $1$ and is identical to the tori of equations (\ref{so2ntor}) and (\ref{so2n+1tor}), respectively; the second has elements of determinant $-1$.

The identity 
\be \label{Bconj} BA(\theta)B^{-1}=A(-\theta) \ee
will become useful below. 

With the maximal tori defined as above, every element of the
orthogonal group can be conjugated to the torus, so each conjugacy
class has a nonempty intersection with the group's maximal torus.  

The counting of conjugacy classes depends on the parity of the order
$m$ of the elements, so we treat the odd and even cases separately.  
\subsection{Odd $m$}
We begin with $N(SO(2n+1),m)$. The block-diagonal matrix diag$(B,
I_{2n-2},-1)$ is an element of $SO(2n+1)$ and equation~(\ref{Bconj}) shows
that conjugation by it takes  $x\in T_{SO(2n+1)}$ to $x'\in
T_{SO(2n+1)}$ where $\theta _1'= -\theta_1$ and the other $\theta_l$
remain the same. Similarly, two elements $x$ and $x'$ of
$T_{SO(2n+1)}$ that differ by $\theta_l'=-\theta_l$ for any
$l=1,\ldots , n$ belong to the same conjugacy class. We therefore
consider only elements of $T_{SO(2n+1)}$ with $\theta _l \in [0,1/2]$
as we did for the symplectic case. As before, we  order the $\theta_l$
to be nondecreasing with $l$.  

For elements of order $m$, we have $\theta _l \in {1\over
  m}(0,1,\ldots , [{m\over 2}])$. So $N(SO(2n+1),m)$ is the number of
weak $([{m\over 2}]+1)$-compositions of $n$.

\begin{formula} \label{oddqinoddso} For any positive integer $n$ and
  any odd integer $m=2k+1$,  
$$ N(SO(2n+1),m)=\binom{n+ \left [ {m\over 2} \right ]}{\left[
    {m\over 2}\right]}.$$
\end{formula}

For $O(2n+1)$, there are two conjugacy classes of maximal tori, i.e.
$T_{SO(2n+1)}$, and $T_{2,\mathrm{odd}}$ in equation~(\ref{o2n+1tor}). However, all
elements of $T_{2,\mathrm{odd}}$ have even order, so none has order
$m=2k+1$. Therefore, the number of conjugacy classes of elements of
odd order in $O(2n+1)$ is the same as that for $SO(2n+1)$, so we get
the following result.
\begin{formula} \label{oddqinoddo} For any positive integer $n$ and
  any odd integer $m=2k+1$,  
$$ N(O(2n+1),m)=\label{nq1} \left ( \begin{array}{c} n+\left [ {m\over 2} \right ]  \\ \left [ {m\over 2} \right ] \end{array} \right ) .$$
\end{formula}
For $O(2n)$, again $T_{2,even}\in T_{O(2n)}$ does not play a role when
$m$ is odd. Also, the block diagonal matrix diag$(B, I_{2n-2})$ is an
element of $O(2n)$, so the results for $O(2n+1)$ and $O(2n)$ are the same.
\begin{formula} \label{oddqineveno} For any positive integer $n$ and
  any odd integer $m=2k+1$,  
$$ N(O(2n),m)=\left ( \begin{array}{c} n+\left [ {m\over 2} \right ] \\ \left [ {m\over 2} \right ] \end{array} \right ) ~.$$
\end{formula}
Things become more subtle for $SO(2n)$: diag$(B,I_{2n-2})$ has
determinant $-1$ so it is not an element of $SO(2n)$. Therefore, it is
no longer the case that if $x,x'\in T_{SO(2n)}$ differ only by
$\theta_i'=-\theta_i$ for some $i$'s then $x$ and $x'$ are necessarily
in the same conjugacy class. However, the block diagonal matrix
diag$(B,B,I_{2n-4})$ \emph{is} in $SO(2n)$, so if  $\theta_l'=-\theta_l$
for an even number of $l$'s, $x$ and $x'$ are in the same conjugacy 
class. 

There are two cases to consider: $\theta_1'=\theta_1=0$ and
$\theta_l\neq 0$  for all $l$. In the first case,
$A(\theta_1)=A(\theta_1')=I_2$, and if $\theta_l'=-\theta_l$ for any
additional $l\geq 2$ (not necessarily an even number of times), then
$x$ and $x'$ are in the same conjugacy class. The number of conjugacy
classes that are represented by elements of $T_{SO(2n)}$ with
$\theta_1=0$ is the number of weak $\left (  \left [ {m\over 2} \right
]+1\right )$-compositions of $n-1$. In the second case $\theta_l\neq 0$
 for all  $l$, the number of classes is the number of weak  $\left [
  {m\over 2} \right ]$-compositions of $n$; since here, flipping the
sign of one $\theta _l$, say $\theta_1'=-\theta_1$ and leaving the
others fixed  lands in  a different conjugacy class, we multiply the
number by two to include all the classes. This leads to the following
theorem. 
\begin{formula} \label{oddqinevenso} For any positive integer $n$ and
  any odd integer $m=2k+1$,  
$$
N(SO(2n),m)=\binom{n+\left [ {m\over 2} \right ]
  -1}{\left [ {m\over 2} \right ]}  +2 \binom{n+\left [ {m\over 2}
    \right ]-1}{ \left [ {m\over 
      2} \right ]-1} = 
\binom{n+\left [ {m\over 2} \right ] -1}{ \left
  [ {m\over 2} \right]} \frac{n+m-1}{n} ~. 
$$
\end{formula}

We now turn to $N(SO(2n+1),m,s)$, where as for the symplectic groups,
$s$ denotes the number of distinct conjugate pairs of eigenvalues of
the elements. For all the orthogonal groups, there are $n$
$\theta_l$'s and $\binom{n-1}{s-1}=\frac sn\binom ns$ ways to
partition them into $s$ nonzero parts. There are 
$\left [{m\over 2}\right ]+1$ possible values for the $\theta _i$. The
same is true for  $O(2n+1)$, and $O(2n)$, yielding the next result.

\begin{formula} For any positive integers $n$ and $s$, and any odd
  integer $m=2k+1$, 
$$  N(SO(2n+1),m,s)=N(O(2n+1),m,s)=N(O(2n),m,s)={s\over n} \binom
ns\binom{\left [ {m\over 2}\right ] +1}{ s}. 
$$
\end{formula}
The above derivation does not apply to $SO(2n)$ because as before,
some classes need to be counted twice due to the absence of $(B,
I_{2n-2})$ in $SO(2n)$. First, we divide the $n$ eigenvalue pairs into
$s$ nonzero parts ($s$-compositions of $n$). In choosing the $s$
eigenvalues out of the $\left [ {m\over 2}\right ] +1$ possibilities,
we differentiate the cases where $\theta_1=0$, which we count once,
from the cases where $\theta_1\neq 0$, which we need to count twice to
account for $\theta_1'=-\theta _1$, $\theta_l'=\theta_l$, $l>1$ which
is in a distinct conjugacy class. We get the following formula.
\begin{formula} \label{soddqinevenso}For any positive integers $n$ and
  $s$ and any odd integer $m=2k+1$,  
\beas  N(SO(2n),m,s)&=&\left ( \begin{array}{c} n-1 \\ s -1 \end{array}
\right ) \left [ \left ( \begin{array}{c} \left [ {m\over 2}\right ]
    \\ s-1 \end{array} \right )+2 \left ( \begin{array}{c} \left [
      {m\over 2}\right ] \\ s \end{array} \right ) \right ] \\ 
&=&{s\over n} \left ( \begin{array}{c} n \\ s  \end{array} \right
)\left ( \begin{array}{c} \left [ {m\over 2}\right ] \\ s \end{array}
\right ){m+1-s\over \left [ {m\over 2}\right ]+1-s}~. 
\eeas
\end{formula}
\subsection{Even m} 
Unlike the case for odd $m$, here we will have to consider $T_2$  in
both $O(2n)$ and $O(2n+1)$. There will also be changes from the odd
$m$ case due to the fact that $\theta_l=1/2$, corresponding to
$A(\theta_l)=-I_2$, can appear. 

For $SO(2n+1)$, we have essentially the same as we did for odd $m$,
i.e. weak $\left ( {m\over 2} +1 \right )$-compositions of $n$.
\begin{formula} \label{evenqinoddso} For any positive integer $n$ and
  any even integer $m=2k$,  
$$ N(SO(2n+1),m)=\left ( \begin{array}{c}
  n+\frac{m}{2}\\ \frac{m}{2} \end{array} \right ).$$
\end{formula}
For $O(2n+1)$, we have to consider conjugacy classes with elements 
whose determinant is $-1$, that is, elements of the second part
$T_{2,\mathrm{odd}}$ of the torus $T_{O(2n+1)}$, not just the elements of determinant $1$ as we did previously. 
But the counting is exactly the
same as in $T_{1,\mathrm{odd}}$, so the next theorem follows.
\begin{formula} \label{evenqinoddo} For any positive integer $n$ and
  any even integer $m=2k$,  
$$ N(O(2n+1),m)=2\left ( \begin{array}{c}
  n+\frac{m}{2}\\ \frac{m}{2}\end{array} \right ) .$$
\end{formula}
Turning to $O(2n)$, we note that elements in $T_{2, even}$ have only
$n-1$ $\theta_l$'s. Other than that, the counting is the same as
before, yielding the next result.
\begin{formula} \label{evenqineveno} For any positive integers $n$ and
  any even integer $m=2k$, 
\beas  N(O(2n),m)&=&\left ( \begin{array}{c}
  n+\frac{m}{2}\\ \frac{m}{2} \end{array} \right )  +\left
( \begin{array}{c} n+\frac{m}{2} -1\\ \frac{m}{2} \end{array} \right )
\\  
&=& \left ( \begin{array}{c} n+\frac{m}{2}
  -1\\ \frac{m}{2} \end{array} \right ){4n+m\over 2n}~. 
\eeas
\end{formula}
For $SO(2n)$, again we need to be careful since $\theta_l'=\pm
\theta_l$  does not always mean $x$ and $x'$ are in the same conjugacy
class. Only when at least one of the $\theta_l$ is $0$ or $1/2$, so
that $A(\theta_l)=\pm I_2$ for that $l$, which commutes with $B$, does
$\theta _l'=\pm\theta_l$ mean  $x$ and $x'$ are in the same conjugacy
class. If no $\theta _l$ is 0 or $1/2$ then if say $\theta
_1'=-\theta_1$ and $\theta_l'=\theta_l$, $l>1$, we have a different
conjugacy class for $x$ and $x'$. The number of conjugacy classes such
that at least one $\theta_l$ is 0 or $1/2$ is the number of weak
$\left ( {m\over 2} +1 \right )$-compositions of $n-1$ (where we have
fixed $\theta_1=0$) plus the number of weak $\left ( {m\over 2} \right
)$-compositions of $n-1$ (where we do not allow $\theta_l=0$ and we
require $\theta_l=1/2$ for some $l$). The number of conjugacy classes
where no $\theta_l$ is 0 or $1/2$ is twice the number of weak $\left (
{m\over 2} -1 \right )$-compositions of $n$. After some algebra we
obtain the next result.

\begin{formula} \label{evenqinevenso} For any positive integer $n$ and
  any even integer $m=2k$, 
$$
N(SO(2n),m)=\left ( \begin{array}{c}
  n+\frac{m}{2}\\ \frac{m}{2} \end{array} \right )  +\left
( \begin{array}{c} n+\frac{m}{2}-2\\  \frac{m}{2}-2 \end{array} \right
) . 
$$
\end{formula}
For $N(SO(2n+1),m,s)$, we have the same calculation as for odd $m$,
and for $N(O(2n+1),m,s)$, we simply double the result to account for
the elements in $T_{2,\mathrm{odd}}$, giving the following formulas.

\begin{formula}For any positive integers $n$ and $s$ and any even integer $m=2k$, 
\beas N(SO(2n+1),m,s)
&=&{s\over n} \left ( \begin{array}{c} n \\ s  \end{array} \right ) \left ( \begin{array}{c}  {m\over 2}+1 \\ s \end{array} \right ); \eeas
\beas
N(O(2n+1),m,s)
&=&{2s\over n} \left ( \begin{array}{c} n \\ s  \end{array} \right ) \left ( \begin{array}{c}  {m\over 2}+1 \\ s \end{array} \right ) .
\eeas
\end{formula}
Next is $O(2n)$, where $T_{2,even}$ has only $n-1$ $\theta_l$'s, so
the contribution from $T_{2,even}$ differs from that from $T_{1,even}$
by replacing $n$ with $n-1$. After some algebra we get the following
theorem. 

\begin{formula}For any positive integers $n$ and $s$ and any even integer $m=2k$, 
$$ N(O(2n),m,s)
={2n-s-1\over n-s}\binom{n-2}{ s -1} \binom{{m\over 2}+1}{ s}. 
$$
\end{formula}
For $SO(2n)$, for each $s$-composition of $n$, the number of conjugacy
classes  of $T_{SO(2n)}$ with $\theta _l\neq 0, 1/2$ for all $l$ is
$\binom {{m\over 2}-1}{s}$ and
the number of conjugacy classes with at least one $\theta _l=0, 1/2$
is the sum of $\binom {{m\over 2}}{s-1}$, which gives the number of conjugacy classes with
$\theta_1=0$, and $\binom { {m\over
    2}-1}{ s-1}$ which gives the number of
conjugacy classes with $\theta_l\neq 0 \; \forall l$ and
$\theta_l=1/2$ for some $l$. As before, we multiply the number for
$\theta_l\neq 0, 1/2$ by 2, and add the rest. After some algebra, we
have our final result.
\begin{formula}For any positive integers $n$ and $s$ and any even
  integer $m=2k$,  
$$
N(SO(2n),m,s)=\left ( \begin{array}{c} n-1 \\ s -1 \end{array} \right ) \left [ \left ( \begin{array}{c}  {m\over 2}+1\\ s \end{array} \right ) + \left ( \begin{array}{c}  {m\over 2}-1\\ s \end{array} \right ) \right ].$$
\end{formula}

\vskip 1cm
\noindent {\bf Acknowledgments}

It is a pleasure to thank Jonathan Pakianathan, Steve Gonek, Ben
Green, and Dragomir Djokovic for discussions. It is also a pleasure to
thank Jonathan Pakianathan for comments on an earlier draft. The work
of the first author was supported in part by US DOE Grant number
DE-FG02-91ER40685, and of the second author in part by NSF Grant
number DMS-1068625.

\end{document}